# Sturm Liouville Problem with Moving Discontinuity Points

Fatma HIRA and Nihat ALTINIŞIK

Department of Mathematics, Science and Arts Faculty, Ondokuz Mayıs University,
55139, Samsun, Turkey

Correspondence should be addressed to Fatma Hıra; fatma.hira@omu.edu.tr

**Abstract:** In this paper, we present a new discontinuous Sturm Liouville problem with symmetrically located discontinuities which are defined depending on a neighborhood of a midpoint of the interval. Also the problem contains an eigenparameter in one of the boundary conditions and has coupled transmission conditions at the discontinuity points. We investigate the properties of the eigenvalues, obtain asymptotic formulas for the eigenvalues and the corresponding eigenfunctions and construct Green's function of this problem.

**Key words:** Sturm-Liouville problems, symmetric and moving discontinuities, eigenvalue and eigenfunction, Green's function

## 1. Introduction

The Sturm-Liouville theory plays an important role in solving many mathematical physics problems [1]. Such research is motivated by theory of heat and mass transfer, vibrating string problems when the string is loaded additionally with point masses. Also, some problems with transmission conditions arise in thermal conduction problems for a thin laminated plate and inhomogenous materials ( see [2,3]; see also wherein references). Heat conduction problems in composite walls have been analysed by several researcher [4-6]. In all of these works, temperature distribution of composite walls involving two or more layers in are investigated spectrally and the results are formulated by the eigenvalues and eigenfunctions of the auxiliary spectral problem. In [7], the presence of infinite number of eigenvalues has been shown and an asymptotic formula has been obtained for the eigenvalues of the spectral problem for the temperature distribution in composite walls. The problem in this work included a two-layer composite wall consisting of different materials, having a common contact surface. When using Fourier's method for heat transfer problems, it becomes necessary to solve a related Sturm Liouville problem in which the separation constant plays the role of a spectral parameter (or eigenparameter). In the case of inhomogeneous materials, the Sturm Liouville equation has variable and not necessarily continuous coefficients so that transmission conditions across the interfaces should be added in the problem. In the present work, we consider a composite medium consisting of two interfaces that are located symmetrically, thus we add to the problem two coupled transmission conditions.

In [8-10] only continuous problems have been investigated. Discontinuous Sturm Liouville problems which contain eigenparameter in one or two boundary conditions and have also transmission conditions (or discontinuity conditions) have been investigated in many works (see [11-20]). Asymptotic formulas for eigenvalues and the corresponding eigenfunctions and computing eigenvalues by sinc method were given in [13-18] and [19,20] for these problems, respectively. Green's function and resolvent operator were constructed and were derived asymptotic approximation formulae for Green function in [14,16].

Throughout this paper we consider the boundary value problem

$$\tau(u) := -u'' + q(x)u = \lambda u, \quad x \in I, \tag{1}$$

with an eigenparameter dependent on one of the boundary conditions

$$B_a(u) := \beta_1 u(a) + \beta_2 u'(a) = 0, \tag{2}$$

$$B_b(u) := \lambda\left(\alpha_1' u(b) - \alpha_2' u'(b)\right) + \left(\alpha_1 u(b) - \alpha_2 u'(b)\right) = 0, \tag{3}$$

and coupled transmission conditions at the points of discontinuity that are $\theta_{-\varepsilon}$ and $\theta_{+\varepsilon}$,

$$T_{-\varepsilon}(u) := u(\theta_{-\varepsilon}+) - \left(\mu_1 u(\theta_{-\varepsilon}-) + \mu_2 u'(\theta_{-\varepsilon}-)\right) = 0, \tag{4}$$

$$T'_{-\varepsilon}(u) := u'(\theta_{-\varepsilon}+) - \left(\mu_1' u(\theta_{-\varepsilon}-) + \mu_2' u'(\theta_{-\varepsilon}-)\right) = 0, \tag{5}$$

$$T_{+\varepsilon}(u) := u(\theta_{+\varepsilon}+) - \left(\eta_1 u(\theta_{+\varepsilon}-) + \eta_2 u'(\theta_{+\varepsilon}-)\right) = 0, \tag{6}$$

$$T'_{+\varepsilon}(u) := u'(\theta_{+\varepsilon}+) - \left(\eta_1' u(\theta_{+\varepsilon}-) + \eta_2' u'(\theta_{+\varepsilon}-)\right) = 0, \tag{7}$$

where $I := [a, \theta_{-\varepsilon}) \cup (\theta_{-\varepsilon}, \theta_{+\varepsilon}) \cup (\theta_{+\varepsilon}, b]$; $0 < \varepsilon < (b-a)/2$, $\theta := (a+b)/2$; $\lambda$ is a spectral parameter; $q(x)$ is a given real valued function which is continuous in $[a, \theta_{-\varepsilon})$, $(\theta_{-\varepsilon}, \theta_{+\varepsilon})$ and $(\theta_{+\varepsilon}, b]$ and has finite limits $q(\theta_{-\varepsilon} \pm 0)$, $q(\theta_{+\varepsilon} \pm 0)$; $\beta_i, \alpha_i, \alpha_i', \mu_i, \mu_i', \eta_i, \eta_i'$ $(i=1,2)$ are real numbers such that $|\beta_1| + |\beta_2| \neq 0$,

$$\rho := \det\begin{pmatrix} \alpha_1' & \alpha_2' \\ \alpha_1 & \alpha_2 \end{pmatrix} > 0, \tag{8}$$

and

$$D_1 := \det\begin{pmatrix} \mu_1 & \mu_2 \\ \mu_1' & \mu_2' \end{pmatrix} > 0, \quad D_2 := \det\begin{pmatrix} \eta_1 & \eta_2 \\ \eta_1' & \eta_2' \end{pmatrix} > 0. \tag{9}$$

Also for convenience we will use the notations $\theta_{\pm\varepsilon}\pm := (\theta \pm \varepsilon) \pm 0$.

We introduce a new discontinuous Sturm Liouville problem with symmetrically located discontinuities which are defined depending on a neighborhood of a midpoint of the interval. $\varepsilon$ is a parameter controlling the change of neighborhood process (it can be called tuning parameter) and by using the change of this $\varepsilon$ parameter it's possible to determine points of discontinuity. That is, two points of discontinuity can be determined in interval $[a,b]$ for each $\varepsilon$ value in the interval $0 < \varepsilon < (b-a)/2$. For example, let $a=-1, b=4$ so that $\theta$ and $\varepsilon$ parameter's interval are $\theta = 7/2$ and $0 < \varepsilon < 5/2$, points of discontinuity are $\theta_{-\varepsilon} = 1$ and $\theta_{+\varepsilon} = 2$ for $\varepsilon = 1/2$, points of discontinuity are $\theta_{-\varepsilon} = -1/2$ and $\theta_{+\varepsilon} = 7/2$ for $\varepsilon = 2$, etc. The main result is that points of discontinuity can be moved by changing $\varepsilon$ parameter, so that they can be called moving discontinuity points. In the special case for our problem when $a=0$, $b=\pi$ and $\varepsilon = d$ $(0 < d < \pi/2)$ is derived in [11,12] these problems do not contain an eigenparameter in the boundary conditions. We extend some classic results of Sturm-Liouville theory to the new symmetric and moving discontinuous case. So that, we define a linear operator $A$ in a suitable Hilbert space $H$ such that the eigenvalues of the problem (1)-(7) coincide with those of $A$, construct a special fundamental system of solutions, obtain the asymptotic formulas for the eigenvalues and the corresponding eigenfunctions that depend on $\varepsilon$ parameter. Finally, we construct Green's function for the problem (1)-(7).

## 2. An Operator Formulation

In this section we will introduce the special inner product in the Hilbert space $L_2(a,b) \oplus \mathbb{C}$ and a symmetric linear operator $A$ defined on this Hilbert space such a way that the problem (1)-(7) can be considered as the eigenvalue problem of this operator.

**Definition 1.** *We define a Hilbert space $H$ of two component vectors by $H := L_2(a,b) \oplus \mathbb{C}$ with the inner product*

$$\langle F,G \rangle_H := \int_a^{\theta_{-\varepsilon}} f(x)\overline{g}(x)dx + \frac{1}{D_1}\int_{\theta_{-\varepsilon}}^{\theta_{+\varepsilon}} f(x)\overline{g}(x)dx + \frac{1}{D_1 D_2}\int_{\theta_{+\varepsilon}}^{b} f(x)\overline{g}(x)dx + \frac{1}{\rho D_1 D_2} f_1 \overline{g_1}$$

(10)

where $F, G \in H$ such that

$$F = (f(x), f_1) = (f(x), R'(f)), \quad G = (g(x), g_1) = (g(x), R'(g)),$$

in which

$$R(f) := \alpha_1 f(b) - \alpha_2 f'(b), \quad R'(f) := \alpha_1' f(b) - \alpha_2' f'(b), \tag{11}$$

$$R(g) := \alpha_1 g(b) - \alpha_2 g'(b), \quad R'(g) := \alpha_1' g(b) - \alpha_2' g'(b). \tag{12}$$

For functions $f(x)$, which is defined on $[a, \theta_{-\varepsilon}) \cup (\theta_{-\varepsilon}, \theta_{+\varepsilon}) \cup (\theta_{+\varepsilon}, b]$ and has finite limits $f(\theta_{-\varepsilon} \pm 0) := \lim_{x \to \theta_{-\varepsilon} \pm 0} f(x)$, $f(\theta_{+\varepsilon} \pm 0) := \lim_{x \to \theta_{+\varepsilon} \pm 0} f(x)$, by $f_{(i)}(x)$ $(i = \overline{1,4})$ we denote the function

$$f_{(1)}(x) := \begin{cases} f(x), & x \in [a, \theta_{-\varepsilon}), \\ f(\theta_{-\varepsilon} -), & x = \theta_{-\varepsilon}, \end{cases} \quad f_{(2)}(x) := \begin{cases} f(\theta_{-\varepsilon} +), & x = \theta_{-\varepsilon}, \\ f(x), & x \in (\theta_{-\varepsilon}, \theta_{+\varepsilon}), \end{cases}$$

$$f_{(3)}(x) := \begin{cases} f(x), & x \in (\theta_{-\varepsilon}, \theta_{+\varepsilon}), \\ f(\theta_{+\varepsilon} -), & x = \theta_{+\varepsilon}, \end{cases} \quad f_{(4)}(x) := \begin{cases} f(\theta_{+\varepsilon} +), & x = \theta_{+\varepsilon}, \\ f(x), & x \in (\theta_{+\varepsilon}, b], \end{cases}$$

which are defined on $I_1 = [a, \theta_{-\varepsilon}], I_2 = [\theta_{-\varepsilon}, \theta_{+\varepsilon}]$ and $I_3 = [\theta_{+\varepsilon}, b]$, respectively.

**Definition 2.** *We define a linear operator $A: D(A) \to H$ by*

$$A \begin{pmatrix} f(x) \\ R'(f) \end{pmatrix} := \begin{pmatrix} \tau(f) \\ -R(f) \end{pmatrix} \quad (13)$$

*where the domain $D(A)$ of the linear operator $A$ is defined as the set of all $F = \begin{pmatrix} f(x) \\ R'(f) \end{pmatrix}$ which is satisfy the conditions $(i)$ $\tau(f) \in L_2(a,b)$, $(ii)$ $f_i(.), f_i'(.)$ are absolutely continuous functions in $I_i$ $(i=1,2,3)$, $(iii)$ $B_a(f) = 0$, $(iv)$ $T_{\pm\varepsilon}(f) = T'_{\pm\varepsilon}(f) = 0$.*

Now we can rewrite the problem (1)-(7) in the operator form as $AF = \lambda F$ where $F = \begin{pmatrix} f(x) \\ R'(f) \end{pmatrix} \in D(A)$.

The eigenvalues and eigenfunctions of the problem (1)-(7) are defined as the eigenvalues and the first components of the corresponding eigenelements of the operator $A$, respectively.

**Lemma 3.** *The operator $A$ in $H$ is symmetric.*

**Proof.** For $F, G \in D(A)$,

$$\langle AF, G \rangle_H = \int_a^{\theta_{-\varepsilon}} \tau(f) \overline{g}(x) dx + \frac{1}{D_1} \int_{\theta_{-\varepsilon}}^{\theta_{+\varepsilon}} \tau(f) \overline{g}(x) dx + \frac{1}{D_1 D_2} \int_{\theta_{+\varepsilon}}^{b} \tau(f) \overline{g}(x) dx - \frac{1}{\rho D_1 D_2} R(f) R'(\overline{g}) \quad (14)$$

By two partial integration we obtain

$$\langle AF, G \rangle_H = \langle F, AG \rangle_H + W(f, \overline{g}; \theta_{-\varepsilon}-) - W(f, \overline{g}; a) + \frac{1}{D_1} W(f, \overline{g}; \theta_{+\varepsilon}-) - \frac{1}{D_1} W(f, \overline{g}; \theta_{-\varepsilon}+) + \frac{1}{D_1 D_2} W(f, \overline{g}; b) - \frac{1}{D_1 D_2} W(f, \overline{g}; \theta_{+\varepsilon}+) - \frac{1}{\rho D_1 D_2} \left( R(f) R'(\overline{g}) - R'(f) R(\overline{g}) \right), \quad (15)$$

where, as usual, by $W(f, g; x)$ we denote the Wronskian of the functions $f$ and $g$:

$$W(f, g; x) = f(x) g'(x) - f'(x) g(x).$$

Since $f$ and $\overline{g}$ satisfy the boundary condition (2), it follows that

$$W(f,\overline{g};a)=0, \tag{16}$$

from the transmission conditions (4)-(7) we get

$$W(f,\overline{g};\theta_{-\varepsilon}+)=D_1 W(f,\overline{g};\theta_{-\varepsilon}-), \tag{17}$$

$$W(f,\overline{g};\theta_{+\varepsilon}+)=D_2 W(f,\overline{g};\theta_{+\varepsilon}-). \tag{18}$$

Further, it is easy to verify that

$$R(f)R'(\overline{g})-R'(f)R(\overline{g})=\rho W(f,\overline{g};b). \tag{19}$$

Finally, substituting (16)-(19) in (15) then we have

$$\langle AF,G \rangle_H = \langle F,AG \rangle_H$$

so $A$ is symmetric. □

**Corollary 4.** *All eigenvalues of the problem (1)-(7) are real.*

## 3. Construction of Fundamental Solutions

**Lemma 5.** *Let the real valued function $q(x)$ be continuous in $[a,b]$ where $f(\lambda), g(\lambda)$ are given entire functions. Then for any $\lambda \in \mathbb{C}$ the equation*

$$-u''+q(x)u=\lambda u, \quad x \in [a,b],$$

*has a unique solution $u = u(x,\lambda)$ such that*

$$u(a)=f(\lambda), u'(a)=g(\lambda) \ \ (or \ u(b)=f(\lambda), u'(b)=g(\lambda)),$$

*and for each $x \in [a,b]$, $u(x,\lambda)$ is an entire function of $\lambda$.*

We will define two solutions

$$\phi_\lambda(x)=\begin{cases}\phi_{1\lambda}(x), & x\in[a,\theta_{-\varepsilon}),\\ \phi_{2\lambda}(x), & x\in(\theta_{-\varepsilon},\theta_{+\varepsilon}),\\ \phi_{3\lambda}(x), & x\in(\theta_{+\varepsilon},b],\end{cases} \quad \chi_\lambda(x)=\begin{cases}\chi_{1\lambda}(x), & x\in[a,\theta_{-\varepsilon}),\\ \chi_{2\lambda}(x), & x\in(\theta_{-\varepsilon},\theta_{+\varepsilon}),\\ \chi_{3\lambda}(x), & x\in(\theta_{+\varepsilon},b],\end{cases}$$

of the equation (1) as follows: Let $\phi_{1\lambda}(x)=\phi_1(x,\lambda)$ be the solution of the equation (1) on $[a,\theta_{-\varepsilon}]$, which satisfies the initial conditions

$$u(a) = \beta_2, \quad u'(a) = -\beta_1 . \tag{20}$$

By virtue of Lemma 5, after defining this solution, we define the solution $\phi_{2\lambda}(x) = \phi_2(x,\lambda)$ of the equation (1) on $[\theta_{-\varepsilon}, \theta_{+\varepsilon}]$ by means of the solution $\phi_{1\lambda}(x)$ by the nonstandard initial conditions

$$u(\theta_{-\varepsilon}+) = \mu_1 \phi_{1\lambda}(\theta_{-\varepsilon}-) + \mu_2 \phi'_{1\lambda}(\theta_{-\varepsilon}-), \quad u'(\theta_{-\varepsilon}+) = \mu'_1 \phi_{1\lambda}(\theta_{-\varepsilon}-) + \mu'_2 \phi'_{1\lambda}(\theta_{-\varepsilon}-) . \tag{21}$$

After defining this solution, we may define the solution $\phi_{3\lambda}(x) = \phi_3(x,\lambda)$ of equation (1) on $[\theta_{+\varepsilon}, b]$ by means of the solution $\phi_{2\lambda}(x)$ by the nonstandard initial conditions

$$u(\theta_{+\varepsilon}+) = \eta_1 \phi_{2\lambda}(\theta_{+\varepsilon}-) + \eta_2 \phi'_{2\lambda}(\theta_{+\varepsilon}-), \quad u'(\theta_{+\varepsilon}+) = \eta'_1 \phi_{2\lambda}(\theta_{+\varepsilon}-) + \eta'_2 \phi'_{2\lambda}(\theta_{+\varepsilon}-) . \tag{22}$$

Hence, $\phi_\lambda(x) = \phi(x,\lambda)$ satisfies the equation (1) on $[a,\theta_{-\varepsilon}) \cup (\theta_{-\varepsilon}, \theta_{+\varepsilon}) \cup (\theta_{+\varepsilon}, b]$, the boundary condition (2) and the transmission conditions (4)-(7).

Analogically, first we define the solution $\chi_{3\lambda}(x) = \chi_3(x,\lambda)$ on $[\theta_{+\varepsilon}, b]$ by the initial conditions

$$u(b) = \alpha'_2 \lambda + \alpha_2, \quad u'(b) = \alpha'_1 \lambda + \alpha_1 . \tag{23}$$

Again, after defining this solution, we define the solution $\chi_{2\lambda}(x) = \chi_2(x,\lambda)$ of the equation (1) on $[\theta_{-\varepsilon}, \theta_{+\varepsilon}]$ by the initial conditions

$$u(\theta_{+\varepsilon}-) = \frac{1}{D_2}\left(\eta'_2 \chi_{3\lambda}(\theta_{+\varepsilon}+) - \eta_2 \chi'_{3\lambda}(\theta_{+\varepsilon}+)\right), \quad u'(\theta_{+\varepsilon}-) = \frac{1}{D_2}\left(-\eta'_1 \chi_{3\lambda}(\theta_{+\varepsilon}+) + \eta_1 \chi'_{3\lambda}(\theta_{+\varepsilon}+)\right). \tag{24}$$

After defining this solution, we define the solution $\chi_{1\lambda}(x) = \chi_1(x,\lambda)$ of the equation (1) on $[a, \theta_{-\varepsilon}]$ by the initial conditions

$$u(\theta_{-\varepsilon}-) = \frac{1}{D_1}\left(\mu'_2 \chi_{2\lambda}(\theta_{-\varepsilon}+) - \mu_2 \chi'_{2\lambda}(\theta_{-\varepsilon}+)\right),$$
$$u'(\theta_{-\varepsilon}-) = \frac{1}{D_1}\left(-\mu'_1 \chi_{2\lambda}(\theta_{-\varepsilon}+) + \mu_1 \chi'_{2\lambda}(\theta_{-\varepsilon}+)\right). \tag{25}$$

Hence, $\chi_\lambda(x) = \chi(x,\lambda)$ satisfies the equation (1) on $[a,\theta_{-\varepsilon}) \cup (\theta_{-\varepsilon}, \theta_{+\varepsilon}) \cup (\theta_{+\varepsilon}, b]$, the boundary condition (3) and the transmission conditions (4)-(7).

Let us consider the Wronskians

$$\omega_i(\lambda) = W(\phi_{i\lambda}, \chi_{i\lambda}; x) = \phi_{i\lambda}(x) \chi'_{i\lambda}(x) - \phi'_{i\lambda}(x) \chi_{i\lambda}(x),$$

which are independent of $x \in I_i$ $(i=1,2,3)$ and entire functions, where $I_i$ $(i=1,2,3)$. After short calculation we see that $D_1 D_2 \omega_1(\lambda) = D_2 \omega_2(\lambda) = \omega_3(\lambda)$. Now we may introduce the characteristic function $\omega(\lambda)$ as

$$\omega(\lambda) := \omega_1(\lambda) = \frac{1}{D_1}\omega_2(\lambda) = \frac{1}{D_1 D_2}\omega_3(\lambda). \tag{26}$$

**Theorem 6.** *The eigenvalues of the problem (1)-(7) are the zeros of the function $\omega(\lambda)$.*

**Proof.** Let $\omega(\lambda_0) = 0$. Then $W(\phi_{1\lambda_0}, \chi_{1\lambda_0}; x) = 0$ and therefore the functions $\phi_{1\lambda_0}(x)$ and $\chi_{1\lambda_0}(x)$ are linearly dependent, i.e. $\chi_{1\lambda_0}(x) = k_1 \phi_{1\lambda_0}(x)$, $x \in [a, \theta_{-\varepsilon}]$ for some $k_1 \neq 0$. From this, it follows that $\chi_{\lambda_0}(x)$ satisfies also the first boundary condition (2), so $\chi_{\lambda_0}(x)$ is an eigenfunction for the eigenvalue $\lambda_0$.

Now let $u_0(x)$ be any eigenfunction correspond to eigenvalue $\lambda_0$, but $\omega(\lambda_0) \neq 0$. Then the pair of the functions $(\phi_1, \chi_1), (\phi_2, \chi_2)$ and $(\phi_3, \chi_3)$ would be linearly independent on $I_i$ $(i=1,2,3)$ respectively. Therefore $u_0(x)$ may be represented as

$$u_0(x) = \begin{cases} c_1 \phi_{1\lambda_0}(x) + c_2 \chi_{1\lambda_0}(x), & x \in [a, \theta_{-\varepsilon}), \\ c_3 \phi_{2\lambda_0}(x) + c_4 \chi_{2\lambda_0}(x), & x \in (\theta_{-\varepsilon}, \theta_{+\varepsilon}), \\ c_5 \phi_{3\lambda_0}(x) + c_6 \chi_{3\lambda_0}(x), & x \in (\theta_{+\varepsilon}, b], \end{cases} \tag{27}$$

where at least one of the constants $c_i$ $(i=\overline{1,6})$ is not zero. Considering the equations

$$B_a(u_0(x)) = 0, \; B_b(u_0(x)) = 0, \; T_{\pm\varepsilon}(u_0(x)) = 0, \; T'_{\pm\varepsilon}(u_0(x)) = 0, \tag{28}$$

as a system of linear equations of the variables $c_i$ $(i=\overline{1,6})$ and taking (21),(22),(24) and (25) into account, it follows that the determinant of this system is

$$\begin{vmatrix} 0 & \omega_1(\lambda_0) & 0 & 0 & 0 & 0 \\ \phi_{2\lambda_0}(\phi_{-\varepsilon}+) & \chi_{2\lambda_0}(\phi_{-\varepsilon}+) & -\phi_{2\lambda_0}(\phi_{-\varepsilon}+) & -\chi_{2\lambda_0}(\phi_{-\varepsilon}+) & 0 & 0 \\ \phi'_{2\lambda_0}(\phi_{-\varepsilon}+) & \chi'_{2\lambda_0}(\phi_{-\varepsilon}+) & -\phi'_{2\lambda_0}(\phi_{-\varepsilon}+) & -\chi'_{2\lambda_0}(\phi_{-\varepsilon}+) & 0 & 0 \\ 0 & 0 & \phi_{3\lambda_0}(\phi_{+\varepsilon}+) & \chi_{3\lambda_0}(\phi_{+\varepsilon}+) & -\phi_{3\lambda_0}(\phi_{+\varepsilon}+) & -\chi_{3\lambda_0}(\phi_{+\varepsilon}+) \\ 0 & 0 & \phi'_{3\lambda_0}(\phi_{+\varepsilon}+) & \chi'_{3\lambda_0}(\phi_{+\varepsilon}+) & -\phi'_{3\lambda_0}(\phi_{+\varepsilon}+) & -\chi'_{3\lambda_0}(\phi_{+\varepsilon}+) \\ 0 & 0 & 0 & 0 & \omega_3(\lambda_0) & 0 \end{vmatrix}$$

$$= -\omega_1(\lambda_0) \omega_2(\lambda_0) \omega_3^2(\lambda_0) \neq 0.$$

Therefore, the system (27) has only the trivial solution $c_i = 0 (i = \overline{1,6})$. Thus we get a contradiction, which completes the proof. □

## 4. Asymptotic Approximate Formulas

It is known [21] that $\phi_\lambda(x)$ be the solution of the equation (1) can be written in the form

$$\phi_{1\lambda}^{(k)}(x) = \beta_2 \left(\cos(s(x-a))\right)^{(k)} - \frac{\beta_1}{s}\left(\sin(s(x-a))\right)^{(k)} +$$
$$\frac{1}{s}\int_a^x \sin(s(x-y))^{(k)} q(y)\phi_{1\lambda}(y)dy, \tag{29}$$

$$\phi_{2\lambda}^{(k)}(x) = \left(\mu_1\phi_{1\lambda}(\theta_{-\varepsilon}-) + \mu_2\phi_{1\lambda}'(\theta_{-\varepsilon}-)\right)\left(\cos(s(x-\theta_{-\varepsilon}))\right)^{(k)} +$$
$$\frac{1}{s}\left(\mu_1'\phi_{1\lambda}(\theta_{-\varepsilon}-) + \mu_2'\phi_{1\lambda}'(\theta_{-\varepsilon}-)\right)\left(\sin(s(x-\theta_{-\varepsilon}))\right)^{(k)} + \tag{30}$$
$$\frac{1}{s}\int_{\theta_{-\varepsilon}}^x \left(\sin(s(x-y))\right)^{(k)} q(y)\phi_{2\lambda}(y)dy,$$

$$\phi_{3\lambda}^{(k)}(x) = \left(\eta_1\phi_{2\lambda}(\theta_{+\varepsilon}-) + \eta_2\phi_{2\lambda}'(\theta_{+\varepsilon}-)\right)\left(\cos(s(x-\theta_{+\varepsilon}))\right)^{(k)} +$$
$$\frac{1}{s}\left(\eta_1'\phi_{2\lambda}(\theta_{+\varepsilon}-) + \eta_2'\phi_{2\lambda}'(\theta_{+\varepsilon}-)\right)\left(\sin(s(x-\theta_{+\varepsilon}))\right)^{(k)} + \tag{31}$$
$$\frac{1}{s}\int_{\theta_{+\varepsilon}}^x \left(\sin(s(x-y))\right)^{(k)} q(y)\phi_{3\lambda}(y)dy,$$

where $s = \sigma + it$, $\lambda = s^2$, $(\cdot)^{(k)} = \dfrac{d^k(\cdot)}{dx^k}$ and for $k = 0,1$.

We also use Titchmarsh's formulas [19] for asymptotic behaviour of $\phi_{i\lambda}(x)(i=\overline{1,3})$ to get the following formulas for the problem (1)-(7) as $\lambda \to \infty$:

$$\phi_{1\lambda}^{(k)}(x) = \beta_2 \left(\cos(s(x-a))\right)^{(k)} + O\left(|s|^{k-1} e^{|t|(x-a)}\right), \tag{32}$$

$$\phi_{2\lambda}^{(k)}(x) = -s\mu_2\beta_2 \sin(s(\theta_{-\varepsilon}-a))\left(\cos(s(x-\theta_{-\varepsilon}))\right)^{(k)} + O\left(|s|^k e^{|t|(x-a)}\right), \tag{33}$$

$$\phi_{3\lambda}^{(k)}(x) = s^2\mu_2\eta_2\beta_2 \sin(s(\theta_{-\varepsilon}-a))\sin(s(\theta_{+\varepsilon}-\theta_{-\varepsilon}))\left(\cos(s(x-\theta_{+\varepsilon}))\right)^{(k)} +$$
$$O\left(|s|^{k+1} e^{|t|(x-a)}\right), \tag{34}$$

for $\beta_2 \neq 0$,

$$\phi_{1\lambda}^{(k)}(x) = -\frac{\beta_1}{s}\left(\sin(s(x-a))\right)^{(k)} + O\left(|s|^{k-2} e^{|t|(x-a)}\right), \tag{35}$$

$$\phi_{2\lambda}^{(k)}(x) = -\mu_2 \beta_1 \cos(s(\theta_{-\varepsilon}-a))\left(\cos(s(x-\theta_{-\varepsilon}))\right)^{(k)} + O\left(|s|^{k-1} e^{|t|(x-a)}\right), \tag{36}$$

$$\phi_{3\lambda}^{(k)}(x) = s\mu_2 \eta_2 \beta_1 \cos(s(\theta_{-\varepsilon}-a))\sin(s(\theta_{+\varepsilon}-\theta_{-\varepsilon}))\left(\cos(s(x-\theta_{+\varepsilon}))\right)^{(k)} + O\left(|s|^{k} e^{|t|(x-a)}\right), \tag{37}$$

for $\beta_2 = 0$.

By substituting the obtained asymptotic formulas for $\phi_{i\lambda}(x)$ $(i=\overline{1,3})$ and the determined solution by (23)-(25) for $\chi_{i\lambda}(x)$ $(i=\overline{1,3})$ into the representation

$$\omega(\lambda) = \frac{1}{D_1 D_2}\left\{(\lambda \alpha_1' + \alpha_1)\phi_{3\lambda}(b) - (\lambda \alpha_2' + \alpha_2)\phi_{3\lambda}'(b)\right\},$$

then $\omega(\lambda)$ has the following asymptotic representations:

**Case 1.** $\beta_2 \neq 0$, $\alpha_2' \neq 0$;

$$\omega(\lambda) = \frac{1}{D_1 D_2} s^5 \beta_2 \alpha_2' \mu_2 \eta_2 \sin(s(\theta_{-\varepsilon}-a))\sin(s(\theta_{+\varepsilon}-\theta_{-\varepsilon}))\sin(s(b-\theta_{+\varepsilon})) + O\left(|s|^4 e^{|t|(b-a)}\right), \tag{38}$$

**Case 2.** $\beta_2 \neq 0$, $\alpha_2' = 0$;

$$\omega(\lambda) = \frac{1}{D_1 D_2} s^4 \beta_2 \alpha_1 \mu_2 \eta_2 \sin(s(\theta_{-\varepsilon}-a))\sin(s(\theta_{+\varepsilon}-\theta_{-\varepsilon}))\cos(s(b-\theta_{+\varepsilon})) + O\left(|s|^3 e^{|t|(b-a)}\right), \tag{39}$$

**Case 3.** $\beta_2 = 0$, $\alpha_2' \neq 0$;

$$\omega(\lambda) = \frac{1}{D_1 D_2} s^4 \beta_1 \alpha_2' \mu_2 \eta_2 \cos(s(\theta_{-\varepsilon}-a))\sin(s(\theta_{+\varepsilon}-\theta_{-\varepsilon}))\sin(s(b-\theta_{+\varepsilon})) + O\left(|s|^3 e^{|t|(b-a)}\right), \tag{40}$$

**Case 4.** $\beta_2 = 0$, $\alpha_2' = 0$;

$$\omega(\lambda) = \frac{1}{D_1 D_2} s^3 \beta_1 \alpha_1 \mu_2 \eta_2 \cos(s(\theta_{-\varepsilon}-a))\sin(s(\theta_{+\varepsilon}-\theta_{-\varepsilon}))\cos(s(b-\theta_{+\varepsilon})) + O\left(|s|^2 e^{|t|(b-a)}\right). \tag{41}$$

**Corollary 7.** *The eigenvalues of the problem (1)-(7) is bounded from below.*

We are now ready to find the asymptotic approximation formulas for the eigenvalues of the problem (1)-(7). Since the eigenvalues coincide with the zeros of the entire functions $\omega(\lambda)$, it follows that they have no finite accumulation point. Morever, all eigenvalues are real and bounded below by Corollaries 4 and 7. Therefore, we may renumber them as $\lambda_0 \leq \lambda_1 \leq \lambda_2 \leq ...$ which counted according to their multiplicity. Below we shall denotes $s_n^2 = \lambda_n$.

**Theorem 8.** *The problem (1)-(7) has an precisely denumerable many real eigenvalues, whose behaviour may be expressed by three sequence $\{\lambda_n'\}, \{\lambda_n''\}$ and $\{\lambda_n'''\}$ with following asymptotics representations for $n \to \infty$:*

**Case 1.** $\beta_2 \neq 0, \alpha_2' \neq 0$;

$$s_n' = \frac{(n-1)\pi}{(\theta_{-\varepsilon} - a)} + O\left(\frac{1}{n}\right), \quad s_n'' = \frac{(n-1)\pi}{(\theta_{+\varepsilon} - \theta_{-\varepsilon})} + O\left(\frac{1}{n}\right), \quad s_n''' = \frac{(n-2)\pi}{(b - \theta_{+\varepsilon})} + O\left(\frac{1}{n}\right), \quad (42)$$

**Case 2.** $\beta_2 \neq 0, \alpha_2' = 0$;

$$s_n' = \frac{(n-1)\pi}{(\theta_{-\varepsilon} - a)} + O\left(\frac{1}{n}\right), \quad s_n'' = \frac{(n-1)\pi}{(\theta_{+\varepsilon} - \theta_{-\varepsilon})} + O\left(\frac{1}{n}\right), \quad s_n''' = \frac{(n-1/2)\pi}{(b - \theta_{+\varepsilon})} + O\left(\frac{1}{n}\right), \quad (43)$$

**Case 3.** $\beta_2 = 0, \alpha_2' \neq 0$;

$$s_n' = \frac{(n-1/2)\pi}{(\theta_{-\varepsilon} - a)} + O\left(\frac{1}{n}\right), \quad s_n'' = \frac{(n-1)\pi}{(\theta_{+\varepsilon} - \theta_{-\varepsilon})} + O\left(\frac{1}{n}\right), \quad s_n''' = \frac{(n-1)\pi}{(b - \theta_{+\varepsilon})} + O\left(\frac{1}{n}\right), \quad (44)$$

**Case 4.** $\beta_2 = 0, \alpha_2' = 0$;

$$s_n' = \frac{(n-1/2)\pi}{(\theta_{-\varepsilon} - a)} + O\left(\frac{1}{n}\right), \quad s_n'' = \frac{(n-1)\pi}{(\theta_{+\varepsilon} - \theta_{-\varepsilon})} + O\left(\frac{1}{n}\right), \quad s_n''' = \frac{(n-1/2)\pi}{(b - \theta_{+\varepsilon})} + O\left(\frac{1}{n}\right). \quad (45)$$

**Proof.** Let us consider Case 1. By applying the well known Rouche's theorem on a sufficiently large contour, it follows that $\omega(\lambda)$ has the same number of zeros inside the contour as the leading term in (38). Hence, if $\lambda_0' \leq \lambda_1' \leq \lambda_2' \leq ...$ are the zeros of $\omega(\lambda)$ and $s_n'^2 = \lambda_n'$ we have

$$s_n' = \frac{(n-1)\pi}{(\theta_{-\varepsilon} - a)} + \delta_n', \quad (46)$$

for sufficiently large $n$, where $\left|\delta_n'\right| \leq \dfrac{\pi}{2(\theta_{-\varepsilon} - a)}$. By putting in (38) we have $\delta_n' = O\left(\dfrac{1}{n}\right)$, which completes the proof for the first formula of the Case 1. The proof for the other cases are similar. □

Then from (32)-(37) (for $k = 0$) and the above theorem, the asymptotic behaviour of the eigenfunctions of the problem (1)-(7) is given by:

**Case 1.** $\beta_2 \neq 0$, $\alpha_2' \neq 0$;

$$\phi_{\lambda_n'}(x) = \begin{cases} \beta_2 \cos\left(\dfrac{(n-1)\pi}{(\theta_{-\varepsilon} - a)}(x - a)\right) + O\left(\dfrac{1}{n}\right), & x \in [a, \theta_{-\varepsilon}), \\ O\left(\dfrac{1}{n}\right), & x \in (\theta_{-\varepsilon}, \theta_{+\varepsilon}), \\ O\left(\dfrac{1}{n}\right), & x \in (\theta_{+\varepsilon}, b], \end{cases}$$

$$\phi_{\lambda_n''}(x) = \begin{cases} \beta_2 \cos\left(\dfrac{(n-1)\pi}{(\theta_{+\varepsilon} - \theta_{-\varepsilon})}(x - a)\right) + O\left(\dfrac{1}{n}\right), & x \in [a, \theta_{-\varepsilon}), \\ -\dfrac{(n-1)\pi}{(\theta_{+\varepsilon} - \theta_{-\varepsilon})}\mu_2\beta_2 \sin\left(\dfrac{(n-1)\pi}{(\theta_{+\varepsilon} - \theta_{-\varepsilon})}(\theta_{-\varepsilon} - a)\right)\cos\left(\dfrac{(n-1)\pi}{(\theta_{+\varepsilon} - \theta_{-\varepsilon})}(x - \theta_{-\varepsilon})\right) + \\ O(1), & x \in (\theta_{-\varepsilon}, \theta_{+\varepsilon}), \\ O\left(\dfrac{1}{n}\right), & x \in (\theta_{+\varepsilon}, b], \end{cases}$$

$$\phi_{\lambda_n'''}(x) = \begin{cases} \beta_2 \cos\left(\dfrac{(n-2)\pi}{(b - \theta_{+\varepsilon})}(x - a)\right) + O\left(\dfrac{1}{n}\right), & x \in [a, \theta_{-\varepsilon}), \\ -\dfrac{(n-2)\pi}{(b - \theta_{+\varepsilon})}\mu_2\beta_2 \sin\left(\dfrac{(n-2)\pi}{(b - \theta_{+\varepsilon})}(\theta_{-\varepsilon} - a)\right)\cos\left(\dfrac{(n-2)\pi}{(b - \theta_{+\varepsilon})}(x - \theta_{-\varepsilon})\right) + \\ O(1), & x \in (\theta_{-\varepsilon}, \theta_{+\varepsilon}), \\ \left(\dfrac{(n-2)\pi}{(b - \theta_{+\varepsilon})}\right)^2 \mu_2\eta_2\beta_2 \sin\left(\dfrac{(n-2)\pi}{(b - \theta_{+\varepsilon})}(\theta_{-\varepsilon} - a)\right)\sin\left(\dfrac{(n-2)\pi}{(b - \theta_{+\varepsilon})}(\theta_{+\varepsilon} - \theta_{-\varepsilon})\right) \times \\ \times \cos\left(\dfrac{(n-2)\pi}{(b - \theta_{+\varepsilon})}(x - \theta_{+\varepsilon})\right) + O(n), & x \in (\theta_{+\varepsilon}, b], \end{cases}$$

**Case 2.** $\beta_2 \neq 0,\ \alpha_2' = 0$;

$$\phi_{\lambda_n}'(x) = \begin{cases} \beta_2 \cos\left(\dfrac{(n-1)\pi}{(\theta_{-\varepsilon}-a)}(x-a)\right) + O\left(\dfrac{1}{n}\right), & x \in [a,\theta_{-\varepsilon}), \\ O\left(\dfrac{1}{n}\right), & x \in (\theta_{-\varepsilon},\theta_{+\varepsilon}), \\ O\left(\dfrac{1}{n}\right), & x \in (\theta_{+\varepsilon},b], \end{cases}$$

$$\phi_{\lambda_n}''(x) = \begin{cases} \beta_2 \cos\left(\dfrac{(n-1)\pi}{(\theta_{+\varepsilon}-\theta_{-\varepsilon})}(x-a)\right) + O\left(\dfrac{1}{n}\right), & x \in [a,\theta_{-\varepsilon}), \\ -\dfrac{(n-1)\pi}{(\theta_{+\varepsilon}-\theta_{-\varepsilon})}\mu_2\beta_2 \sin\left(\dfrac{(n-1)\pi}{(\theta_{+\varepsilon}-\theta_{-\varepsilon})}(\theta_{-\varepsilon}-a)\right)\cos\left(\dfrac{(n-1)\pi}{(\theta_{+\varepsilon}-\theta_{-\varepsilon})}(x-\theta_{-\varepsilon})\right) + \\ O(1), & x \in (\theta_{-\varepsilon},\theta_{+\varepsilon}), \\ O\left(\dfrac{1}{n}\right), & x \in (\theta_{+\varepsilon},b], \end{cases}$$

$$\phi_{\lambda_n}''(x) = \begin{cases} \beta_2 \cos\left(\dfrac{(n-1/2)\pi}{(b-\theta_{+\varepsilon})}(x-a)\right) + O\left(\dfrac{1}{n}\right), & x \in [a,\theta_{-\varepsilon}), \\ -\dfrac{(n-1/2)\pi}{(b-\theta_{+\varepsilon})}\mu_2\beta_2 \sin\left(\dfrac{(n-1/2)\pi}{(b-\theta_{+\varepsilon})}(\theta_{-\varepsilon}-a)\right)\cos\left(\dfrac{(n-1/2)\pi}{(b-\theta_{+\varepsilon})}(x-\theta_{-\varepsilon})\right) + \\ O(1), & x \in (\theta_{-\varepsilon},\theta_{+\varepsilon}), \\ \left(\dfrac{(n-1/2)\pi}{(b-\theta_{+\varepsilon})}\right)^2 \mu_2\eta_2\beta_2 \sin\left(\dfrac{(n-1/2)\pi}{(b-\theta_{+\varepsilon})}(\theta_{-\varepsilon}-a)\right)\sin\left(\dfrac{(n-1/2)\pi}{(b-\theta_{+\varepsilon})}(\theta_{+\varepsilon}-\theta_{-\varepsilon})\right) \times \\ \times \cos\left(\dfrac{(n-1/2)\pi}{(b-\theta_{+\varepsilon})}(x-\theta_{+\varepsilon})\right) + O(n), & x \in (\theta_{+\varepsilon},b], \end{cases}$$

**Case 3.** $\beta_2 = 0,\ \alpha_2' \neq 0$;

$$\phi_{\lambda_n}'(x) = \begin{cases} -\dfrac{\beta_1(\theta_{-\varepsilon}-a)}{(n-1/2)\pi}\sin\left(\dfrac{(n-1/2)\pi}{(\theta_{-\varepsilon}-a)}(x-a)\right) + O\left(\dfrac{1}{n^2}\right), & x \in [a,\theta_{-\varepsilon}), \\ O\left(\dfrac{1}{n}\right), & x \in (\theta_{-\varepsilon},\theta_{+\varepsilon}), \\ O(1), & x \in (b,\theta_{+\varepsilon}], \end{cases}$$

$$\phi_{\lambda_n^*}(x) = \begin{cases} -\dfrac{\beta_1(\theta_{+\varepsilon}-\theta_{-\varepsilon})}{(n-1)\pi}\sin\left(\dfrac{(n-1)\pi}{(\theta_{+\varepsilon}-\theta_{-\varepsilon})}(x-a)\right)+O\left(\dfrac{1}{n^2}\right), & x\in[a,\theta_{-\varepsilon}), \\ -\mu_2\beta_1\cos\left(\dfrac{(n-1)\pi}{(\theta_{+\varepsilon}-\theta_{-\varepsilon})}(\theta_{-\varepsilon}-a)\right)\cos\left(\dfrac{(n-1)\pi}{(\theta_{+\varepsilon}-\theta_{-\varepsilon})}(x-\theta_{-\varepsilon})\right)+O\left(\dfrac{1}{n}\right), & x\in(\theta_{-\varepsilon},\theta_{+\varepsilon}), \\ O(1), & x\in(b,\theta_{+\varepsilon}], \end{cases}$$

$$\phi_{\lambda_n^*}(x) = \begin{cases} -\dfrac{\beta_1(b-\theta_{+\varepsilon})}{(n-1)\pi}\sin\left(\dfrac{(n-1)\pi}{(b-\theta_{+\varepsilon})}(x-a)\right)+O\left(\dfrac{1}{n^2}\right), & x\in[a,\theta_{-\varepsilon}), \\ -\mu_2\beta_1\cos\left(\dfrac{(n-1)\pi}{(b-\theta_{+\varepsilon})}(\theta_{-\varepsilon}-a)\right)\cos\left(\dfrac{(n-1)\pi}{(b-\theta_{+\varepsilon})}(x-\theta_{-\varepsilon})\right)+O\left(\dfrac{1}{n}\right), & x\in(\theta_{-\varepsilon},\theta_{+\varepsilon}), \\ \dfrac{(n-1)\pi}{(b-\theta_{+\varepsilon})}\mu_2\beta_1\eta_2\cos\left(\dfrac{(n-1)\pi}{(b-\theta_{+\varepsilon})}(\theta_{-\varepsilon}-a)\right)\sin\left(\dfrac{(n-1)\pi}{(b-\theta_{+\varepsilon})}(\theta_{+\varepsilon}-\theta_{-\varepsilon})\right)\times \\ \times\cos\left(\dfrac{(n-1)\pi}{(b-\theta_{+\varepsilon})}(x-\theta_{-\varepsilon})\right)+O(1), & x\in(b,\theta_{+\varepsilon}], \end{cases}$$

**Case 4.** $\beta_2 = 0$, $\alpha_2' = 0$;

$$\phi_{\lambda_n^*}(x) = \begin{cases} -\dfrac{\beta_1(\theta_{-\varepsilon}-a)}{(n-1/2)\pi}\sin\left(\dfrac{(n-1/2)\pi}{(\theta_{-\varepsilon}-a)}(x-a)\right)+O\left(\dfrac{1}{n^2}\right), & x\in[a,\theta_{-\varepsilon}), \\ O\left(\dfrac{1}{n}\right), & x\in(\theta_{-\varepsilon},\theta_{+\varepsilon}), \\ O(1), & x\in(b,\theta_{+\varepsilon}], \end{cases}$$

$$\phi_{\lambda_n^*}(x) = \begin{cases} -\dfrac{\beta_1(\theta_{+\varepsilon}-\theta_{-\varepsilon})}{(n-1)\pi}\sin\left(\dfrac{(n-1)\pi}{(\theta_{+\varepsilon}-\theta_{-\varepsilon})}(x-a)\right)+O\left(\dfrac{1}{n^2}\right), & x\in[a,\theta_{-\varepsilon}), \\ -\mu_2\beta_1\cos\left(\dfrac{(n-1)\pi}{(\theta_{+\varepsilon}-\theta_{-\varepsilon})}(\theta_{-\varepsilon}-a)\right)\cos\left(\dfrac{(n-1)\pi}{(\theta_{+\varepsilon}-\theta_{-\varepsilon})}(x-\theta_{-\varepsilon})\right)+O\left(\dfrac{1}{n}\right), & x\in(\theta_{-\varepsilon},\theta_{+\varepsilon}), \\ O(1), & x\in(b,\theta_{+\varepsilon}], \end{cases}$$

$$\phi_{\lambda_n^*}(x) = \begin{cases} -\dfrac{\beta_1(b-\theta_{+\varepsilon})}{(n-1/2)\pi}\sin\left(\dfrac{(n-1/2)\pi}{(b-\theta_{+\varepsilon})}(x-a)\right)+O\left(\dfrac{1}{n^2}\right), & x\in[a,\theta_{-\varepsilon}), \\[2mm] -\mu_2\beta_1\cos\left(\dfrac{(n-1/2)\pi}{(b-\theta_{+\varepsilon})}(\theta_{-\varepsilon}-a)\right)\cos\left(\dfrac{(n-1/2)\pi}{(b-\theta_{+\varepsilon})}(x-\theta_{-\varepsilon})\right)+O\left(\dfrac{1}{n}\right), & x\in(\theta_{-\varepsilon},\theta_{+\varepsilon}), \\[2mm] \dfrac{(n-1/2)\pi}{(b-\theta_{+\varepsilon})}\mu_2\beta_1\eta_2\cos\left(\dfrac{(n-1/2)\pi}{(b-\theta_{+\varepsilon})}(\theta_{-\varepsilon}-a)\right)\sin\left(\dfrac{(n-1/2)\pi}{(b-\theta_{+\varepsilon})}(\theta_{+\varepsilon}-\theta_{-\varepsilon})\right)\times \\[2mm] \times\cos\left(\dfrac{(n-1/2)\pi}{(b-\theta_{+\varepsilon})}(x-\theta_{-\varepsilon})\right)+O(1), & x\in(b,\theta_{+\varepsilon}], \end{cases}$$

All these asymptotic approximations hold uniformly for $x\in I$.

## 5. Green's Function

Now let $\lambda\in\mathbb{C}$ not be an eigenvalue of $A$ and consider the inhomogenous problem for $F=\begin{pmatrix}f(x)\\f_1\end{pmatrix}\in H$, $U=\begin{pmatrix}u(x)\\R'(u)\end{pmatrix}\in D(A)$,

$$(\lambda I - A)U = F, \ x\in I \tag{47}$$

and I is the identity operator. Since

$$(\lambda I - A)U = \lambda\begin{pmatrix}u(x)\\R'(u)\end{pmatrix}-\begin{pmatrix}\tau(u)\\-R(u)\end{pmatrix}=\begin{pmatrix}f(x)\\f_1\end{pmatrix}, \tag{48}$$

then we have

$$(\lambda-\tau)u(x) = f(x), \ x\in I \tag{49}$$

$$\lambda R'(u)+R(u) = f_1, \tag{50}$$

Now we can represent the general solution of homogeneous differential equation (1), appropriate to equation (49) in the following form:

$$u(x,\lambda) = \begin{cases} c_1\phi_{1\lambda}(x)+c_2\chi_{1\lambda}(x), & x\in[a,\theta_{-\varepsilon}), \\ c_3\phi_{2\lambda}(x)+c_4\chi_{2\lambda}(x), & x\in(\theta_{-\varepsilon},\theta_{+\varepsilon}), \\ c_5\phi_{3\lambda}(x)+c_6\chi_{3\lambda}(x), & x\in(\theta_{+\varepsilon},b], \end{cases}$$

in which $c_i, i=\overline{1,6}$ are arbitrary constants. By applying the method of variation of the constants, we shall search the general solution of the non-homogeneous linear differential equation (49) in the following form:

$$u(x,\lambda) = \begin{cases} c_1(x,\lambda)\phi_{1\lambda}(x) + c_2(x,\lambda)\chi_{1\lambda}(x), & x \in [a, \theta_{-\varepsilon}), \\ c_3(x,\lambda)\phi_{2\lambda}(x) + c_4(x,\lambda)\chi_{2\lambda}(x), & x \in (\theta_{-\varepsilon}, \theta_{+\varepsilon}), \\ c_5(x,\lambda)\phi_{3\lambda}(x) + c_6(x,\lambda)\chi_{3\lambda}(x), & x \in (\theta_{+\varepsilon}, b], \end{cases} \quad (51)$$

where the functions $c_i(x,\lambda)$ $i = \overline{1,6}$ satisfy the linear system of equations

$$\begin{cases} c_1'(x,\lambda)\phi_{1\lambda}(x) + c_2'(x,\lambda)\chi_{1\lambda}(x) = 0, \\ c_1'(x,\lambda)\phi'_{1\lambda}(x) + c_2'(x,\lambda)\chi'_{1\lambda}(x) = f(x), \end{cases} \text{ for } x \in [a, \theta_{-\varepsilon}), \quad (52)$$

$$\begin{cases} c_3'(x,\lambda)\phi_{2\lambda}(x) + c_4'(x,\lambda)\chi_{2\lambda}(x) = 0, \\ c_3'(x,\lambda)\phi'_{2\lambda}(x) + c_4'(x,\lambda)\chi'_{2\lambda}(x) = f(x), \end{cases}, \text{ for } x \in (\theta_{-\varepsilon}, \theta_{+\varepsilon}), \quad (53)$$

$$\begin{cases} c_5'(x,\lambda)\phi_{3\lambda}(x) + c_6'(x,\lambda)\chi_{2\lambda}(x) = 0, \\ c_5'(x,\lambda)\phi'_{3\lambda}(x) + c_6'(x,\lambda)\chi'_{3\lambda}(x) = f(x), \end{cases}, \text{ for } x \in (\theta_{-\varepsilon}, b]. \quad (54)$$

Since $\lambda$ is not an eigenvalue and $\omega_i(\lambda) \neq 0$ $(i = 1,2,3)$ each of the linear systems in (52)-(54) have a unique solution which leads

$$c_1(x,\lambda) = \frac{1}{\omega_1(\lambda)} \int_x^{\theta_{-\varepsilon}} \chi_{1\lambda}(y) f(y) dy + c_1(\lambda), \quad c_2(x,\lambda) = \frac{1}{\omega_1(\lambda)} \int_a^x \phi_{1\lambda}(y) f(y) dy + c_2(\lambda),$$

$$c_3(x,\lambda) = \frac{1}{\omega_2(\lambda)} \int_x^{\theta_{+\varepsilon}} \chi_{2\lambda}(y) f(y) dy + c_3(\lambda), \quad c_4(x,\lambda) = \frac{1}{\omega_2(\lambda)} \int_{\theta_{-\varepsilon}}^x \phi_{2\lambda}(y) f(y) dy + c_4(\lambda),$$

$$c_5(x,\lambda) = \frac{1}{\omega_3(\lambda)} \int_x^b \chi_{3\lambda}(y) f(y) dy + c_5(\lambda), \quad c_6(x,\lambda) = \frac{1}{\omega_3(\lambda)} \int_{\theta_{+\varepsilon}}^x \phi_{3\lambda}(y) f(y) dy + c_6(\lambda),$$

$$(55)$$

where $c_i(\lambda)$ $i = \overline{1,6}$ are arbitrary constants. Substituting (55) into (51), we obtain the solution of (49),

$$u(x,\lambda) = \begin{cases} \dfrac{\phi_{1\lambda}(x)}{\omega_1(\lambda)} \displaystyle\int_x^{\theta_{-\varepsilon}} \chi_{1\lambda}(y) f(y) dy + \dfrac{\chi_{1\lambda}(x)}{\omega_1(\lambda)} \displaystyle\int_a^x \phi_{1\lambda}(y) f(y) dy + \\ c_1(\lambda)\phi_{1\lambda}(x) + c_2(\lambda)\chi_{1\lambda}(x), \qquad x \in [a, \theta_{-\varepsilon}), \\[1em] \dfrac{\phi_{2\lambda}(x)}{\omega_2(\lambda)} \displaystyle\int_x^{\theta_{+\varepsilon}} \chi_{2\lambda}(y) f(y) dy + \dfrac{\chi_{2\lambda}(x)}{\omega_2(\lambda)} \displaystyle\int_{\theta_{-\varepsilon}}^x \phi_{2\lambda}(y) f(y) dy + \\ c_3(\lambda)\phi_{2\lambda}(x) + c_4(\lambda)\chi_{2\lambda}(x), \qquad x \in (\theta_{-\varepsilon}, \theta_{+\varepsilon}), \\[1em] \dfrac{\phi_{3\lambda}(x)}{\omega_3(\lambda)} \displaystyle\int_x^{b} \chi_{3\lambda}(y) f(y) dy + \dfrac{\chi_{3\lambda}(x)}{\omega_3(\lambda)} \displaystyle\int_{\theta_{+\varepsilon}}^x \phi_{3\lambda}(y) f(y) dy + \\ c_5(\lambda)\phi_{3\lambda}(x) + c_6(\lambda)\chi_{3\lambda}(x), \qquad x \in (\theta_{+\varepsilon}, b], \end{cases} \qquad (56)$$

Then from boundary conditions (2), (50) and transmission conditions (4)-(7) we get

$$c_1(\lambda) = \frac{1}{\omega_2(\lambda)} \int_{\theta_{-\varepsilon}}^{\theta_{+\varepsilon}} \chi_{2\lambda}(y) f(y) dy + \frac{1}{\omega_3(\lambda)} \int_{\theta_{+\varepsilon}}^{b} \chi_{3\lambda}(y) f(y) dy + \frac{f_1}{\omega_3(\lambda)}, \quad c_2(\lambda) = 0,$$

$$c_3(\lambda) = \frac{1}{\omega_3(\lambda)} \int_{\theta_{+\varepsilon}}^{b} \chi_{3\lambda}(y) f(y) dy + \frac{f_1}{\omega_3(\lambda)}, \quad c_4(\lambda) = \frac{1}{\omega_1(\lambda)} \int_{a}^{\theta_{-\varepsilon}} \phi_{1\lambda}(y) f(y) dy,$$

$$c_5(\lambda) = \frac{f_1}{\omega_3(\lambda)}, \quad c_6(\lambda) = \frac{1}{\omega_1(\lambda)} \int_{a}^{\theta_{-\varepsilon}} \phi_{1\lambda}(y) f(y) dy + \frac{1}{\omega_2(\lambda)} \int_{\theta_{-\varepsilon}}^{\theta_{+\varepsilon}} \phi_{2\lambda}(y) f(y) dy.$$

(57)

Substituting (57) and (26) into (56), then (56) can be written as

$$u(x,\lambda) = \begin{cases} \dfrac{\phi_{1\lambda}(x)}{\omega(\lambda)} \int\limits_{x}^{\theta_{-\varepsilon}} \chi_{1\lambda}(y) f(y) dy + \dfrac{\chi_{1\lambda}(x)}{\omega(\lambda)} \int\limits_{a}^{x} \phi_{1\lambda}(y) f(y) dy + \\ \dfrac{\phi_{1\lambda}(x)}{D_1 \omega(\lambda)} \int\limits_{\theta_{-\varepsilon}}^{\theta_{+\varepsilon}} \chi_{2\lambda}(y) f(y) dy + \dfrac{\phi_{1\lambda}(x)}{D_1 D_2 \omega(\lambda)} \int\limits_{\theta_{+\varepsilon}}^{b} \chi_{3\lambda}(y) f(y) dy + \\ \dfrac{f_1}{D_1 D_2 \omega(\lambda)} \phi_{1\lambda}(x), \qquad\qquad\qquad\qquad x \in [a, \theta_{-\varepsilon}), \\[4pt]
\dfrac{\phi_{2\lambda}(x)}{D_1 \omega(\lambda)} \int\limits_{x}^{\theta_{+\varepsilon}} \chi_{2\lambda}(y) f(y) dy + \dfrac{\chi_{2\lambda}(x)}{D_1 \omega(\lambda)} \int\limits_{\theta_{-\varepsilon}}^{x} \phi_{2\lambda}(y) f(y) dy + \\ \dfrac{\chi_{2\lambda}(x)}{\omega(\lambda)} \int\limits_{a}^{\theta_{-\varepsilon}} \phi_{1\lambda}(y) f(y) dy + \dfrac{\phi_{2\lambda}(x)}{D_1 D_2 \omega(\lambda)} \int\limits_{\theta_{+\varepsilon}}^{b} \chi_{3\lambda}(y) f(y) dy + \\ \dfrac{f_1}{D_1 D_2 \omega(\lambda)} \phi_{2\lambda}(x), \qquad\qquad\qquad\qquad x \in (\theta_{-\varepsilon}, \theta_{+\varepsilon}), \\[4pt]
\dfrac{\phi_{3\lambda}(x)}{D_1 D_2 \omega(\lambda)} \int\limits_{x}^{b} \chi_{3\lambda}(y) f(y) dy + \dfrac{\chi_{3\lambda}(x)}{D_1 D_2 \omega(\lambda)} \int\limits_{\theta_{+\varepsilon}}^{x} \phi_{3\lambda}(y) f(y) dy + \\ \dfrac{\chi_{3\lambda}(x)}{\omega(\lambda)} \int\limits_{a}^{\theta_{-\varepsilon}} \phi_{1\lambda}(y) f(y) dy + \dfrac{\chi_{3\lambda}(x)}{D_1 \omega(\lambda)} \int\limits_{\theta_{-\varepsilon}}^{\theta_{+\varepsilon}} \phi_{2\lambda}(y) f(y) dy + \\ \dfrac{f_1}{D_1 D_2 \omega(\lambda)} \phi_{3\lambda}(x), \qquad\qquad\qquad\qquad x \in (\theta_{+\varepsilon}, b]. \end{cases} \quad (58)$$

Then (58) can be rewritten in the form

$$u(x,\lambda) = \int\limits_{a}^{\theta_{-\varepsilon}} G(x,y;\lambda) f(y) dy + \dfrac{1}{D_1} \int\limits_{\theta_{-\varepsilon}}^{\theta_{+\varepsilon}} G(x,y;\lambda) f(y) dy + \\ \dfrac{1}{D_1 D_2} \int\limits_{\theta_{+\varepsilon}}^{b} G(x,y;\lambda) f(y) dy + \dfrac{f_1 \phi_\lambda(x)}{D_1 D_2 \omega(\lambda)}, \quad (59)$$

where

$$G(x,y;\lambda) = \begin{cases} \dfrac{\phi_\lambda(y) \chi_\lambda(x)}{\omega(\lambda)}, & a \le y \le x \le b,\ x \ne \theta_{-\varepsilon}, \theta_{+\varepsilon};\ y \ne \theta_{-\varepsilon}, \theta_{+\varepsilon}, \\[6pt] \dfrac{\phi_\lambda(x) \chi_\lambda(y)}{\omega(\lambda)}, & a \le x \le y \le b,\ x \ne \theta_{-\varepsilon}, \theta_{+\varepsilon};\ y \ne \theta_{-\varepsilon}, \theta_{+\varepsilon}, \end{cases} \quad (60)$$

is Green's function of the problem (1)-(7).

Hence, we have

$$U = (\lambda I - A)^{-1} F =$$

$$\left( \int_a^{\theta_{-\varepsilon}} G(x,y;\lambda) f(y) dy + \frac{1}{D_1} \int_{\theta_{-\varepsilon}}^{\theta_{+\varepsilon}} G(x,y;\lambda) f(y) dy + \frac{1}{D_1 D_2} \int_{\theta_{+\varepsilon}}^b G(x,y;\lambda) f(y) dy + \frac{f_1 \phi_\lambda(x)}{D_1 D_2 \omega(\lambda)} \right)$$

$$R'(u)$$

(61)

the resolvent of the problem (1)-(7).

**Conflict of interests.** The authors declare that there is no conflict of interests regarding the publication of this paper.

**Acknowledgement.** This work was supported by the Research Formulation of Ondokuz Mayıs University under the project number PYO.FEN.1904.11.010.